\documentclass[a4paper,conference]{IEEEtran}
\IEEEoverridecommandlockouts

\usepackage{cite,color,url}
\usepackage{amsmath,amssymb,amsfonts,hyperref,multirow}
\usepackage{booktabs}
\usepackage{algorithmic,tikz}
\usepackage{subcaption,amsthm,longtable,lscape,amsmath}
\usepackage{graphicx,longtable}
\usepackage{textcomp}
\usepackage{xcolor}
\usetikzlibrary{fit,calc}
\def\BibTeX{{\rm B\kern-.05em{\sc i\kern-.025em b}\kern-.08em
    T\kern-.1667em\lower.7ex\hbox{E}\kern-.125emX}}
\begin{document}

\title{On Solving the Knapsack Problem with Conflicts}
\author{\IEEEauthorblockN{Roberto Montemanni}
\IEEEauthorblockA{\textit{Department of Sciences and Methods for Engineering} \\
\textit{University of Modena and Reggio Emilia}\\
Reggio Emilia, Italy \\
roberto.montemanni@unimore.it}
\and
\IEEEauthorblockN{Derek H. Smith}
\IEEEauthorblockA{\textit{Faculty of Computing, Engineering and Science} \\
\textit{University of South Wales}\\
Pontypridd, Wales, UK \\
derek.smith@southwales.ac.uk}}

\maketitle
\begin{abstract}
A variant of the well-known Knapsack Problem is studied in this paper, where pairs of items are conflicting, and cannot be selected at the same time. This configures a set of hard constraints. The problem, which can be used to model real applications, looks for a selection of items such that the total profit is maximized, the capacity of the container is respected, and no conflict is violated. 

In this paper, we consider a previously known mixed integer linear
program representing the problem and we solve it with the open-source solver CP-SAT, part of the Google OR-Tools computational suite.

An experimental campaign on the instances available from the literature and adopted in the last decade, indicate that the approach we propose achieves results comparable with, and often better than, those of state-of-the-art solvers, notwithstanding its intrinsic conceptual and implementation simplicity. 
\end{abstract}


\section{Introduction}
The classic Knapsack Problem \cite{mar90} consists of selecting objects from a given set, where each object is associated with a profit and a weight, in such a way that the total weight of the objects selected is below a given threshold and the total profit is maximized. In this work a generalization of the problem, known as the Knapsack Problem with Conflicts (KPC), where pairs $(i, j)$ of conflicting objects are given, and for each pair at most one object can be part of the final solution.

The KPC arises as a subproblem in several algorithms in the Operations Research domain. For example, it is solved within branch-and-price methods for the Bin Packing Problem \cite{wei19} and the Bin Packing Problem with Conflicts \cite{bet17}. However, the KPC can be used to model in mathematical terms some real-world applications. Less abstract applications can be found in the scheduling domain, where tasks that cannot be executed in parallel due to conflicting resources can be mapped into conflicting items.

The KPC was originally introduced  in \cite{yam02}, where a local-search heuristic and a branch-and-bound method based on Lagrangian relaxation are discussed. In \cite{hif07} the authors proposed some preprocessing techniques and a different branch-and bound algorithm.  Heuristic methods for the KPC were discussed in \cite{hif06} and \cite{ake11}, the latter method being a matheuristic algorithm. Metaheuristic approaches were finally presented in \cite{hif12} (scatter-search) and \cite{hif14} (rounding heuristic). Finally, a branch-and bound exact method based on binary branching was proposed in \cite{bet17}, while another branch-and-bound schema based on an $n$-ary branching was developed and tested in \cite{con21}.

In this paper, a mixed integer linear programming model for the KPC is considered and solved by the open-source solver CP-SAT, which is part of the Google OR-Tools \cite{cpsat} optimization suite. Successful application of this solver on optimization problems with  characteristics similar to the problem under investigation, motivated our study \cite{md23}, \cite{cor}, \cite{rm25}. An experimental campaign on the  benchmark instances previously proposed in the recent literature is also presented and discussed.

The overall organization of the paper can be summarized as follows. The Knapack Problem with Conflicts is formally defined in Section \ref{desc}. Section \ref{model} discusses a mixed integer linear programming model to represent the problem. In Section \ref{exp} the model is tested experimentally on the benchmark instances adopted in the literature in the last decade. The approach we propose is compared with recent state-of-the-art methods from the literature. Conclusions are finally drawn in Section \ref{conc}.

\section{Problem Description}\label{desc}
The KPC can be formally defined as follows. Let $G=(V,E)$ be a graph where each vertex $i \in V$ is associated with an item, and each edge $\{i,j\} \in E$ models a conflict between the items $i$ and $j$ of $V$. A profit $p_i$ and a weight $w_i$ are provided for each item $i \in V$, and a capacity $c$ of the knapsack is finally provided. The aim of the problem is to select a subset $S$ of the items of $V$ such that the total weight of the items of $S$ is not greater than $c$,  no conflict is violated (i.e. $\forall i,j \in S, \{i,j\} \notin E$), and the total profit of the selected items is maximized.  

\begin{figure*}[h!]
{
\begin{center}
{
\begin{tikzpicture}[node distance={2cm}, main/.style = {draw, circle}]
			\node at (-2,0.4) {\textcolor{red}{$p_i$}, \textcolor{blue}{$w_i$}};
			\node at (-2,-0.25) {\textcolor{blue}{$c=20$}};
			\node[main,minimum size=1.15cm] (0) {\textcolor{red}{6}, \textcolor{blue}{7}};
			\node[main,minimum size=1.15cm] (1) [right of=0] {\textcolor{red}{9}, \textcolor{blue}{9}};
			\node[main,minimum size=1.15cm] (2) [right of =1] {\textcolor{red}{9}, \textcolor{blue}{4}};
			\node[main,minimum size=1.15cm] (3) [below of=0] {\textcolor{red}{3}, \textcolor{blue}{3}};
			\node[main,minimum size= 1.15cm] (4) [right of =3] {\textcolor{red}{7}, \textcolor{blue}{6}};
			\node[main,minimum size=1.15cm] (5) [right of=4] {\textcolor{red}{2}, \textcolor{blue}{1}};	
			\draw [thick, line width=1.2,-] (0) to node {} (1) ;
			\draw [thick, line width=1.2,-] (0) to node {} (4) ;
			\draw [thick, line width=1.2,-] (0) to node {} (1) ;
			\draw [thick, line width=1.2,-] (2) to node {} (2) ;
			\draw [thick, line width=1.2,-] (3) to node {} (2) ;
			\draw [thick, line width=1.2,-] (4) to node {} (2) ;
\end{tikzpicture}
\hspace{2cm}
\begin{tikzpicture}[node distance={2cm}, main/.style = {draw, circle}]
			\node at (6,0.4) {\textcolor{red}{$\sum p_i =21$}};
			\node at (6,-0.25) {\textcolor{blue}{$\sum w_i =19$}};
			\node[main,minimum size=1.15cm] (1) [right of=0] {\textcolor{red}{9}, \textcolor{blue}{9}};
			\node[main,minimum size=1.15cm] (3) [below of=0] {\textcolor{red}{3}, \textcolor{blue}{3}};
			\node[main,minimum size= 1.15cm] (4) [right of =3] {\textcolor{red}{7}, \textcolor{blue}{6}};
			\node[main,minimum size=1.15cm] (5) [right of=4] {\textcolor{red}{2}, \textcolor{blue}{1}};	
\end{tikzpicture}
}
	\caption{On the left an example of a KPC instance is presented, where for each item the profit is indicated in red and the weight in blue and edges represent conflicts. The total capacity $c$ is 20. On the right, an optimal solution with total profit 21 and weight 19 is represented. Observe that no conflict is violated by the given solution.}
	\label{figu}
\end{center}
}
\end{figure*}
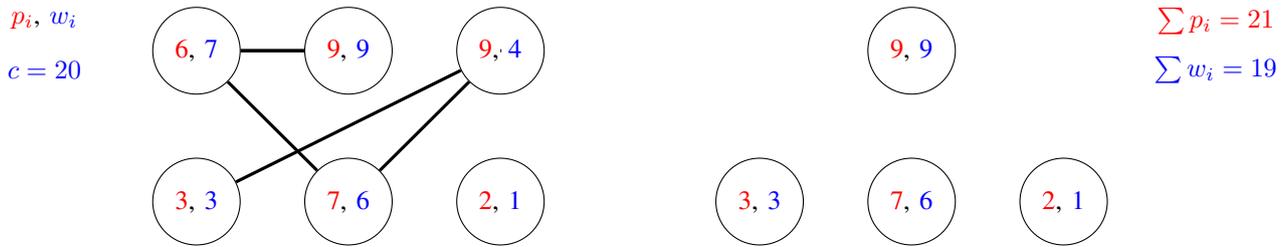

A small example of a KPC instance and a related optimal solution are depicted in Figure \ref{figu}.

\section{A Mixed Integer Linear Programming Model}\label{model}
In this section, a model for the KPC, previously discussed in \cite{bet17} and in other works, is presented. A variable $x_{i}$ takes value 1 if the item $ i\in V$ is selected, 0 otherwise. The resulting model is as follows.

\begin{align} 
  \max \ \ &   \sum_{i \in V} p_i x_i& \label{1}\\ 
s.t. \ \ &	\sum_{i \in V}  w_i x_{i} \le c& \label{2}\\
& x_i + x_j \le 1& \{i,j\} \in E \label{3}\\
& x_{i} \in \{0,1\} & j \in V \label{4}
\end{align}

The objective function (\ref{1}) maximizes the profit of the items selected.
Constraints (\ref{2}) impose that the total weight of the selected items cannot exceed the given capacity $c$. 
Inequalities (\ref{3}) model the conflicts, imposing that at most one of two conflicting items can be selected.
The domains of the $x$ variables are finally specified in constraints (\ref{4}).

\section{Computational Experiments} \label{exp}
In Section \ref{ben} we describe the benchmark instances previously introduced in the literature, and used for the present study. In Section \ref{res} the approach we propose is compared with the other methods available in the literature.

\subsection{Benchmark Instances}\label{ben}
The first benchmark set adopted for the experiments was proposed in \cite{bet17}. The instances are derived from those originally introduced in \cite{fal96}  for the
 bin packing problem. They consist of eight classes, each composed by
10 instances; in the first four classes the items have weights that follow a uniform distribution in the range $[20, 100]$ and
the capacity $c$ is 150. The number $n$ of
items is 120, 250, 500, and 1,000, respectively. The last
four classes have weights sampled with a uniform distribution
in the range $[250,500]$ and a capacity of 1000. 
The instances have 60, 120, 349, and 501 items, respectively.
A random conflict graph is generated for each instance, with density values in
the range from 0.1 to 0.9, making a total of 90 instances for each class. A profit
is associated with each item. Profits are either uniformly distributed in the interval $[1, 100]$ ($R$ instances), or defined as $p_i = w_i + 10$ for each item $i$ ($C$ instances). Three copies of each instances are obtained by considering the original capacity ($C1$, $R1$), or a capacity obtained by multiplying the original one by 3 ($C3$, $R3$) or 10 ($C10$, $R10$).
Overall, 4320 instances compose the first dataset.

The second group of instances -- again first introduced in \cite{bet17} -- considers very sparse
conflict graphs, with densities 0.001, 0.002, 0.005 and 0.01. Capacities are either 1000 or 2000 and the number of items is 500 or 1000. Ten instances with random profits and 10 with correlated profits were generated for each combination of parameters, for a total of 480 instances.

We refer the interested reader to \cite{bet17} for a comprehensive description of the instances.

\subsection{Experimental Results} \label{res}
The model discussed in Section \ref{model} has been solved with the Google OR-Tools CP-SAT solver \cite{cpsat} version 9.12. The experiments have been run on a computer equipped with an Intel Core i7 12700F CPU. The experiments for the methods previously appeared in \cite{con21} and against which we compare, were run on an a machine equipped with an Intel Xeon CPU E5-2690, which according to \url{http://gene.disi.unitn.it/test/cpu_list.php} is approximately 10\% slower. Moreover, in \cite{con21} the experiments were run on a single core, while we let the solver use all the potentialities of the processor. These two factors give an advantage to the approach we propose, although this will not change the outcome of the general experiments.  

The methods involved in the comparison are:
\begin{itemize}
\item CFS: a combinatorial branch-and-bound algorithm, presented in \cite{con21};
\item BCM: a branch-and-bound algorithm based on binary branching and a strong upper bounding procedure, presented in \cite{bet17};
\item ILP: three Integer Linear Programming models discussed in \cite{con21} and solved with ILOG CPLEX 12.8 \cite{cplex}. For each dataset the results of the best of the three models are presented;
\item CP-SAT: the mixed integer linear program presented in Section \ref{model} solved with Google OR-Tools CP-SAT solver 9.12 \cite{cpsat}.
\end{itemize}

All the methods considered were run for a maximum of 600 seconds on each instance, and for each group of instances the number of proven optimal solution and the average time to find such optimal solutions, are provided. For the method CP-SAT, the average optimality gap is also reported for the second benchmark set, in order to measure the quality of the solutions also when optimality could not be proven (this information was not available for the other methods). Results in bold highlight improvements of CP-SAT over the state-of-the-art.

The results on the first benchmark set are presented in Tables \ref{t11} and \ref{t12}, with data are aggregated by Class and Type and Class and Density, respectively.

\begin{table}[]
{
\caption{Results on the first benchmark set (aggregated by Class and Type)}\label{t11}
\begin{tabular}{lrrrrrrrrr}
\toprule
\multicolumn{2}{c}{Instances}&  \multicolumn{2}{c}{CFS \cite{con21}}          & \multicolumn{2}{c}{BCM \cite{bet17}}        &\multicolumn{2}{c}{ILP \cite{con21}}       & \multicolumn{2}{c}{CP-SAT}      \\
Class\!\!\!\! & T   & Opt & Sec  & Opt  & Sec   & Opt & Sec   & Opt & Sec  \\
\cmidrule(lr){1-2}\cmidrule(lr){3-4}\cmidrule(lr){5-6}\cmidrule(lr){7-8}\cmidrule(lr){9-10}
C1  & 1 & 90 & 0.0  & 90 & 0.0   & 90 & 0.2   & 90 & 0.1  \\
    & 2 & 90 & 0.0  & 90 & 0.0   & 90 & 1.1   & 90 & 0.3  \\
    & 3 & 90 & 0.0  & 90 & 0.0   & 90 & 8.2   & 90 & 1.7  \\
    & 4 & 90 & 0.0  & 90 & 0.0   & 90 & 24.2  & 90 & 14.4 \\
    & 5 & 90 & 0.0  & 90 & 0.0   & 90 & 0.0   & 90 & 0.0  \\
    & 6 & 90 & 0.0  & 90 & 0.0   & 90 & 0.1   & 90 & 0.1  \\
    & 7 & 90 & 0.0  & 90 & 0.0   & 90 & 0.5   & 90 & 0.2  \\
    & 8 & 90 & 0.0  & 90 & 0.0   & 90 & 3.6   & 90 & 0.7  \\
        \midrule
C3  & 1 & 90 & 0.0  & 90 & 0.0   & 90 & 1.5   & 90 & 0.0  \\
    & 2 & 90 & 0.0  & 90 & 0.1   & 90 & 25.8  & 90 & 0.1  \\
    & 3 & 90 & 0.1  & 90 & 1.3   & 54 & 162.8 & 90 & 0.2  \\
    & 4 & 90 & 1.6  & 90 & 27.3  & 21 & 141.9 & 90 & 1.9  \\
    & 5 & 90 & 0.0  & 90 & 0.0   & 90 & 0.2   & 90 & 0.0  \\
    & 6 & 90 & 0.0  & 90 & 0.0   & 90 & 2.0   & 90 & 0.0  \\
    & 7 & 90 & 0.0  & 90 & 0.1   & 90 & 46.5  & 90 & 0.1  \\
    & 8 & 90 & 0.0  & 90 & 0.6   & 59 & 35.3  & 90 & 0.2  \\
        \midrule
C10 & 1 & 90 & 0.1  & 90 & 1.6   & 90 & 3.5   & 90 & 0.1  \\
    & 2 & 90 & 25.2 & 73 & 31.9  & 68 & 126.2 & 90 & \textbf{0.3}  \\
    & 3 & 61 & 15.9 & 50 & 18.2  & 22 & 166.0 & \textbf{90} & \textbf{1.6}  \\
    & 4 & 50 & 47.2 & 40 & 108.8 & 1  & 575.1 & \textbf{90} & \textbf{14.5} \\
    & 5 & 90 & 0.0  & 90 & 0.0   & 90 & 0.2   & 90 & 0.0  \\
    & 6 & 90 & 0.5  & 90 & 6.8   & 90 & 5.3   & 90 & \textbf{0.1}  \\
    & 7 & 86 & 35.9 & 70 & 24.5  & 65 & 143.1 & \textbf{90} & \textbf{0.2}  \\
    & 8 & 60 & 7.3  & 49 & 17.4  & 20 & 156.4 & \textbf{90} & \textbf{0.7}  \\
        \midrule
R1  & 1 & 90 & 0.0  & 90 & 0.0   & 90 & 0.1   & 90 & 0.1  \\
    & 2 & 90 & 0.0  & 90 & 0.0   & 90 & 0.8   & 90 & 0.2  \\
    & 3 & 90 & 0.0  & 90 & 0.0   & 90 & 4.8   & 90 & 1.0  \\
    & 4 & 90 & 0.0  & 90 & 0.1   & 90 & 10.1  & 90 & 9.0  \\
    & 5 & 90 & 0.0  & 90 & 0.0   & 90 & 0.0   & 90 & 0.0  \\
    & 6 & 90 & 0.0  & 90 & 0.0   & 90 & 0.1   & 90 & 0.1  \\
    & 7 & 90 & 0.0  & 90 & 0.0   & 90 & 0.4   & 90 & 0.2  \\
    & 8 & 90 & 0.0  & 90 & 0.1   & 90 & 2.7   & 90 & 0.8  \\
        \midrule
R3  & 1 & 90 & 0.0  & 90 & 0.0   & 90 & 0.4   & 90 & 0.0  \\
    & 2 & 90 & 0.0  & 90 & 0.0   & 90 & 5.0   & 90 & 0.1  \\
    & 3 & 90 & 0.0  & 90 & 0.2   & 90 & 55.1  & 90 & 0.2  \\
    & 4 & 90 & 0.1  & 90 & 2.3   & 50 & 127.2 & 90 & 1.9  \\
    & 5 & 90 & 0.0  & 90 & 0.0   & 90 & 0.1   & 90 & 0.0  \\
    & 6 & 90 & 0.0  & 90 & 0.0   & 90 & 0.5   & 90 & 0.0  \\
    & 7 & 90 & 0.0  & 90 & 0.0   & 90 & 5.0   & 90 & 0.1  \\
    & 8 & 90 & 0.0  & 90 & 0.2   & 90 & 64.7  & 90 & 0.2  \\
    \midrule
R10 & 1 & 90 & 0.0  & 90 & 0.1   & 90 & 1.6   & 90 & 0.1  \\
    & 2 & 90 & 0.8  & 90 & 9.1   & 87 & 107.4 & 90 & \textbf{0.2}  \\
    & 3 & 89 & 49.9 & 69 & 57.0  & 33 & 100.0 & \textbf{90} & \textbf{1.0}  \\
    & 4 & 51 & 23.2 & 40 & 25.0  & 8  & 333.6 & \textbf{90} & \textbf{9.0}  \\
    & 5 & 90 & 0.0  & 90 & 0.0   & 90 & 0.1   & 90 & 0.0  \\
    & 6 & 90 & 0.0  & 90 & 0.2   & 90 & 1.5   & 90 & 0.1  \\
    & 7 & 90 & 1.5  & 90 & 17.7  & 80 & 91.7  & 90 & \textbf{0.2}  \\
    & 8 & 79 & 19.5 & 69 & 43.2  & 30 & 77.0  & \textbf{90} & \textbf{0.8} \\
    \midrule
   \multicolumn{2}{c}{Average}&    86.8&	4.8&	84.6&	8.2	&76.2	&54.6	&\textbf{90.0}	&	\textbf{1.3}\\
    \bottomrule
\end{tabular}
}
\end{table}

\begin{table}[]
\caption{Results on the first benchmark set (aggregated by Class and Density)}\label{t12}
\begin{tabular}{lrrrrrrrrr}
\toprule
\multicolumn{2}{c}{Instances}&  \multicolumn{2}{c}{CFS \cite{con21}}          & \multicolumn{2}{c}{BCM \cite{bet17}}        &\multicolumn{2}{c}{ILP \cite{con21}}       & \multicolumn{2}{c}{CP-SAT}      \\
Class\!\!\!\! & D      & Opt & Sec  & Opt  & Sec   & Opt & Sec   & Opt & Sec  \\
\cmidrule(lr){1-2}\cmidrule(lr){3-4}\cmidrule(lr){5-6}\cmidrule(lr){7-8}\cmidrule(lr){9-10}
C1  & 0.1 & 80 & 0.0  & 80 & 0.0  & 80 & 0.1   & 80 & 0.2 \\
    & 0.2 & 80 & 0.0  & 80 & 0.0  & 80 & 0.2   & 80 & 0.4 \\
    & 0.3 & 80 & 0.0  & 80 & 0.0  & 80 & 0.4   & 80 & 0.7 \\
    & 0.4 & 80 & 0.0  & 80 & 0.0  & 80 & 0.6   & 80 & 1.2 \\
    & 0.5 & 80 & 0.0  & 80 & 0.0  & 80 & 1.3   & 80 & 1.5 \\
    & 0.6 & 80 & 0.0  & 80 & 0.0  & 80 & 2.7   & 80 & 3.2 \\
    & 0.7 & 80 & 0.0  & 80 & 0.0  & 80 & 5.5   & 80 & 3.3 \\
    & 0.8 & 80 & 0.0  & 80 & 0.0  & 80 & 14.4  & 80 & 5.2 \\
    & 0.9 & 80 & 0.0  & 80 & 0.0  & 80 & 17.5  & 80 & 4.0 \\
   \midrule
C3  & 0.1 & 80 & 0.2  & 80 & 0.3  & 77 & 11.2  & 80 & 0.2 \\
    & 0.2 & 80 & 0.1  & 80 & 2.4  & 79 & 26.7  & 80 & 0.4 \\
    & 0.3 & 80 & 0.4  & 80 & 6.6  & 72 & 27.7  & 80 & 0.7 \\
    & 0.4 & 80 & 0.5  & 80 & 9.1  & 66 & 41.3  & 80 & 1.2 \\
    & 0.5 & 80 & 0.5  & 80 & 9.6  & 52 & 74.6  & 80 & 1.6 \\
    & 0.6 & 80 & 0.2  & 80 & 3.9  & 50 & 41.9  & 80 & 3.3 \\
    & 0.7 & 80 & 0.1  & 80 & 1.0  & 50 & 11.0  & 80 & 3.3 \\
    & 0.8 & 80 & 0.0  & 80 & 0.2  & 66 & 68.4  & 80 & 5.2 \\
    & 0.9 & 80 & 0.0  & 80 & 0.0  & 72 & 27.2  & 80 & 3.9 \\
      \midrule
C10 & 0.1 & 47 & 41.4 & 33 & 37.3 & 47 & 53.3  & \textbf{80} & \textbf{0.1} \\
    & 0.2 & 50 & 79.0 & 30 & 4.4  & 30 & 4.5   & \textbf{80} & \textbf{0.1} \\
    & 0.3 & 50 & 1.6  & 50 & 64.7 & 30 & 3.7   & \textbf{80} & \textbf{0.1} \\
    & 0.4 & 70 & 11.5 & 50 & 3.7  & 48 & 169.0 & \textbf{80} & \textbf{0.2} \\
    & 0.5 & 80 & 28.8 & 69 & 23.6 & 50 & 106.4 & 80 & \textbf{0.4} \\
    & 0.6 & 80 & 1.2  & 80 & 52.5 & 50 & 38.2  & 80 & \textbf{0.4} \\
    & 0.7 & 80 & 0.1  & 80 & 3.7  & 50 & 12.4  & 80 & 0.5 \\
    & 0.8 & 80 & 0.0  & 80 & 0.3  & 70 & 75.8  & 80 & 0.5 \\
    & 0.9 & 80 & 0.0  & 80 & 0.0  & 71 & 28.7  & 80 & 0.6 \\
      \midrule
R1  & 0.1 & 80 & 0.0  & 80 & 0.0  & 80 & 0.1   & 80 & 0.3 \\
    & 0.2 & 80 & 0.0  & 80 & 0.0  & 80 & 0.2   & 80 & 0.4 \\
    & 0.3 & 80 & 0.0  & 80 & 0.0  & 80 & 0.3   & 80 & 0.8 \\
    & 0.4 & 80 & 0.0  & 80 & 0.0  & 80 & 0.6   & 80 & 1.2 \\
    & 0.5 & 80 & 0.0  & 80 & 0.0  & 80 & 0.9   & 80 & 1.5 \\
    & 0.6 & 80 & 0.0  & 80 & 0.0  & 80 & 1.6   & 80 & 1.8 \\
    & 0.7 & 80 & 0.0  & 80 & 0.0  & 80 & 4.8   & 80 & 2.1 \\
    & 0.8 & 80 & 0.0  & 80 & 0.0  & 80 & 5.3   & 80 & 2.6 \\
    & 0.9 & 80 & 0.0  & 80 & 0.0  & 80 & 7.2   & 80 & 2.2 \\
       \midrule
R3  & 0.1 & 80 & 0.0  & 80 & 0.1  & 80 & 0.1   & 80 & 0.3 \\
    & 0.2 & 80 & 0.0  & 80 & 0.1  & 80 & 1.0   & 80 & 0.4 \\
    & 0.3 & 80 & 0.0  & 80 & 0.4  & 80 & 10.9  & 80 & 0.8 \\
    & 0.4 & 80 & 0.0  & 80 & 0.7  & 78 & 35.6  & 80 & 1.2 \\
    & 0.5 & 80 & 0.0  & 80 & 0.8  & 70 & 23.3  & 80 & 1.6 \\
    & 0.6 & 80 & 0.0  & 80 & 0.6  & 70 & 46.3  & 80 & 1.8 \\
    & 0.7 & 80 & 0.0  & 80 & 0.3  & 70 & 53.4  & 80 & 2.1 \\
    & 0.8 & 80 & 0.0  & 80 & 0.1  & 74 & 45.2  & 80 & 2.5 \\
    & 0.9 & 80 & 0.0  & 80 & 0.0  & 78 & 31.2  & 80 & 2.3 \\
      \midrule
R10 & 0.1 & 71 & 11.2 & 69 & 59.9 & 72 & 16.0  & \textbf{80} & \textbf{0.1} \\
    & 0.2 & 59 & 47.9 & 50 & 39.5 & 43 & 85.1  & \textbf{80} & \textbf{0.1} \\
    & 0.3 & 70 & 42.8 & 50 & 3.9  & 44 & 139.2 & \textbf{80} & \textbf{0.1} \\
    & 0.4 & 70 & 2.0  & 70 & 38.9 & 50 & 53.8  & \textbf{80} & \textbf{0.2} \\
    & 0.5 & 80 & 7.2  & 70 & 3.9  & 50 & 51.6  & 80 & \textbf{0.4} \\
    & 0.6 & 80 & 0.4  & 80 & 11.6 & 50 & 16.2  & 80 & 0.4 \\
    & 0.7 & 80 & 0.1  & 80 & 1.3  & 54 & 38.5  & 80 & 0.5 \\
    & 0.8 & 80 & 0.0  & 80 & 0.2  & 70 & 40.4  & 80 & 0.5 \\
    & 0.9 & 79 & 0.0  & 79 & 0.0  & 75 & 44.6  & \textbf{80} & \textbf{0.6}\\
    \midrule
   \multicolumn{2}{c}{Average}&77.9	&6.0	&76.2	&8.4	&69.5	&34.6	&\textbf{80.0}	&	\textbf{1.2}\\
    \bottomrule
\end{tabular}
\end{table}

The results on the first benchmark set clearly highlight the superiority of the approach we propose, based on CP-SAT. The dominance does not appear to be affected by the details of the instances, and suggests that for densities in any non-extreme range, the CP-SAT solver performs better than the others. All the instances are solved to optimality and with average computation times orders of magnitude smaller than those of the other solvers. 
 
The results on the second benchmark set are presented in Tables \ref{t13} for the correlated instances and in Table \ref{t14} for the random instances, with data are aggregated by Items/Capacity and Density.

\begin{table*}[]
\begin{center}
\caption{Results on the set benchmark set of \emph{correlated} instances (aggregated by number of Items/Capacity and Density)}\label{t13}
\begin{tabular}{lrrrrrrrrrr}
\toprule
\multicolumn{2}{c}{Instances}&  \multicolumn{2}{c}{CFS \cite{con21}}          & \multicolumn{2}{c}{BCM \cite{bet17}}        &\multicolumn{2}{c}{ILP \cite{con21}}       & \multicolumn{3}{c}{CP-SAT}      \\
Items/Cap\!\!\!\! & Density      & Opt & Sec  & Opt  & Sec   & Opt & Sec   & Opt & Gap \% & Sec  \\
\cmidrule(lr){1-2}\cmidrule(lr){3-4}\cmidrule(lr){5-6}\cmidrule(lr){7-8}\cmidrule(lr){9-11}
500/1000  & 0.001 & 10  & 0.0   & 10  & 0.0   & 10  & 0.0   & 10  & 0.00 & 0.1     \\
          & 0.002 & 10  & 0.0   & 10  & 0.6   & 10  & 0.0   & 10  & 0.00 & 0.1     \\
          & 0.005 & 10  & 0.2   & 10  & 6.7   & 10  & 0.0   & 10  & 0.00 & 13.1    \\
          & 0.01  & 10  & 0.8   & 9   & 103.3 & 10  & 0.0   & 9   & 0.02 & 30.4    \\
          & 0.02  & 10  & 56.7  & 1   & 272.7 & 10  & 0.3   & 10  & 0.00 & 24.8    \\
          & 0.05  & 1   & 165.8 & 0   & -     & 10  & 90.6  & 9   & 0.09 & 204.0   \\
          \midrule
500/2000  & 0.001 & 10  & 4.2   & 10  & 0.4   & 10  & 0.0   & 10  & 0.00 & 0.1     \\
          & 0.002 & 10  & 0.1   & 10  & 5.2   & 10  & 0.0   & 10  & 0.00 & 0.2     \\
          & 0.005 & 10  & 7.3   & 8   & 199.6 & 10  & 0.0   & 9   & 0.01 & 100.9   \\
          & 0.01  & 7   & 49.8  & 0   & -     & 10  & 0.0   & 6   & 0.04 & 5.3     \\
          & 0.02  & 0   & -     & 0   & -     & 9   & 6.1   & 9   & 0.01 & 12.2    \\
          & 0.05  & 0   & -     & 0   & -     & 0   & -     & 0   & 2.77 & -       \\
          \midrule
1000/1000 & 0.001 & 10  & 0.1   & 10  & 5.4   & 10  & 0.0   & 10  & 0.00 & 0.3     \\
          & 0.002 & 10  & 0.2   & 10  & 11.1  & 10  & 0.0   & 10  & 0.00 & 0.5     \\
          & 0.005 & 10  & 5.9   & 5   & 379.9 & 10  & 0.0   & 7   & 0.05 & 22.5 \\
          & 0.01  & 7   & 163.9 & 0   & -     & 10  & 0.1   & 2   & 0.14 & 84.6 \\
          & 0.02  & 0   & -     & 0   & -     & 7   & 2.4   & 5   & 0.11 & 2.9  \\
          & 0.05  & 0   & -     & 0   & -     & 0   & -     & 0   & 4.72 & -       \\
          \midrule
1000/2000 & 0.001 & 10  & 3.1   & 9   & 84.7  & 10  & 0.0   & 10  & 0.00 & 8.7     \\
          & 0.002 & 10  & 45.8  & 7   & 210.3 & 10  & 0.0   & 8   & 0.02 & 37.0 \\
          & 0.005 & 7   & 182.0 & 0   & -     & 10  & 0.0   & 6   & 0.03 & 159.8 \\
          & 0.01  & 4   & 0.0   & 0   & -     & 9   & 0.1   & 6   & 0.04 & 100.4 \\
          & 0.02  & 0   & -     & 0   & -     & 5   & 193.8 & 1   & 0.82 & 388.1 \\
          & 0.05  & 0   & -     & 0   & -     & 0   & -     & 0   & 8.21 & -       \\
          \midrule
\multicolumn{2}{c}{Average}       & 6.5 & 38.1  & 4.5 & 98.5  & 8.3 & 14.0  & 7.0 & 0.71 & 56.9  \\
   \bottomrule
\end{tabular}
\end{center}
\end{table*}

\begin{table*}[]
\begin{center}
\caption{Results on the set benchmark set of \emph{random} instances (aggregated by number of Items/Capacity and Density)}\label{t14}
\begin{tabular}{lrrrrrrrrrr}
\toprule
\multicolumn{2}{c}{Instances}&  \multicolumn{2}{c}{CFS \cite{con21}}          & \multicolumn{2}{c}{BCM \cite{bet17}}        &\multicolumn{2}{c}{ILP \cite{con21}}       & \multicolumn{3}{c}{CP-SAT}      \\
Items/Cap\!\!\!\! & Density      & Opt & Sec  & Opt  & Sec   & Opt & Sec   & Opt & Gap \% & Sec  \\
\cmidrule(lr){1-2}\cmidrule(lr){3-4}\cmidrule(lr){5-6}\cmidrule(lr){7-8}\cmidrule(lr){9-11}
500/1000  & 0.001 & 10  & 0.0   & 10  & 0.0   & 10  & 0.0   & 10  & 0.00  & 0.1     \\
          & 0.002 & 10  & 0.0   & 10  & 0.1   & 10  & 0.0   & 10  & 0.00  & 0.1     \\
          & 0.005 & 10  & 0.0   & 10  & 0.4   & 10  & 0.0   & 10  & 0.00  & 0.2     \\
          & 0.01  & 10  & 0.1   & 10  & 2.1   & 10  & 0.0   & 10  & 0.00  & 0.1     \\
          & 0.02  & 10  & 1.2   & 10  & 32.8  & 10  & 0.0   & 10  & 0.00  & 0.2     \\
          & 0.05  & 9   & 132.7 & 3   & 133.2 & 10  & 1.2   & 10  & 0.00  & 1.2     \\
          \midrule
500/2000  & 0.001 & 10  & 0.0   & 10  & 0.1   & 10  & 0.0   & 10  & 0.00  & 0.1     \\
          & 0.002 & 10  & 0.0   & 10  & 0.3   & 10  & 0.0   & 10  & 0.00  & 0.2     \\
          & 0.005 & 10  & 0.1   & 10  & 2.4   & 10  & 0.0   & 10  & 0.00  & 0.2     \\
          & 0.01  & 10  & 10.4  & 9   & 190.7 & 10  & 0.0   & 10  & 0.00  & 0.1     \\
          & 0.02  & 3   & 116.5 & 1   & 39.6  & 10  & 0.1   & 10  & 0.00  & 0.2     \\
          & 0.05  & 0   & -     & 0   & -     & 10  & 81.2  & 6   & 1.61  & 229.0 \\
          \midrule
1000/1000 & 0.001 & 10  & 0.0   & 10  & 0.4   & 10  & 0.0   & 10  & 0.0   & 0.2     \\
          & 0.002 & 10  & 0.0   & 10  & 1.6   & 10  & 0.0   & 10  & 0.0   & 0.4     \\
          & 0.005 & 10  & 0.1   & 10  & 16.8  & 10  & 0.0   & 10  & 0.0   & 0.3     \\
          & 0.01  & 10  & 15.0  & 8   & 152.6 & 10  & 0.1   & 10  & 0.0   & 0.2     \\
          & 0.02  & 4   & 125.7 & 1   & 468.8 & 10  & 0.6   & 10  & 0.0   & \textbf{0.4}     \\
          & 0.05  & 0   & -     & 0   & -     & 9   & 255.0 & 3   & 2.9   & 139.3   \\
          \midrule
1000/2000 & 0.001 & 10  & 0.0   & 10  & 2.3   & 10  & 0.0   & 10  & 0.00  & 0.3     \\
          & 0.002 & 10  & 0.0   & 10  & 20.3  & 10  & 0.0   & 10  & 0.00  & 0.6     \\
          & 0.005 & 9   & 69.5  & 2   & 189.9 & 10  & 0.0   & 10  & 0.00  & 0.3     \\
          & 0.01  & 1   & 565.8 & 0   & -     & 10  & 0.1   & 10  & 0.00  & 0.3     \\
          & 0.02  & 0   & -     & 0   & -     & 10  & 2.1   & 10  & 0.00  & 10.2    \\
          & 0.05  & 0   & -     & 0   & -     & 0   & -     & 0   & 12.45 & -       \\
          \midrule
\multicolumn{2}{c}{Average}       & 7.3 & 51.9  & 6.4 & 66.0  & 9.5 & 14.8  & 9.1 & 0.71  & 16.7\\
          \bottomrule
\end{tabular}
\end{center}
\end{table*}

The results on both correlated and random instances of the second benchmark set suggest that on these problems, characterized by extreme low densities of the conflict graph, the method we propose remains competitive but is slightly inferior to the ILP approach. In particular, it appears that the most critical values of densities are those around 0.02 and 0.05 for all methods. Larger instances with higher capacities are, on the other hand, the most difficult.

In conclusion, CP-SAT performs better than the other methods on most of the instances, but for extreme cases such as densities such as 0.02 and 0.05, and large instances, its superiority vanishes in favour of ILP.

\section{Conclusions} \label{conc}
A  formulation based on mixed integer linear formulation for the Knapsack Problem with Conflicts has been considered and solved via the open-source solver CP-SAT, part of the Google OR-Tools computational suite. 

The experimental results indicate that the approach we propose has performances comparable with, and often better than, those of the state-of-the-art solvers available in the literature, notwithstanding the minimal implementation effort required for our solution.

\section*{Acknowledgments}
The work was partially supported by the Google Cloud Research Credits Program.

\bibliographystyle{IEEEtran}
\bibliography{mybibfile}

\end{document}